\documentclass{article}

\usepackage[utf8]{inputenc}
\usepackage{amsmath}
\usepackage{amssymb}
\usepackage{amsthm}
\usepackage{bbm}
\usepackage{tikz}

\newtheorem{theorem}{Theorem}

\newcommand{\Prob}{\operatorname{\mathbb{P}}}
\newcommand{\E}{\operatorname{E}}
\newcommand{\cov}{\operatorname{Cov}}
\newcommand{\Var}{\operatorname{Var}}
\newcommand{\aut}{\operatorname{Aut}}

\title{The Probability of Random Trees Being Isomorphic}
\author{Christoffer Olsson}
\date{ }

\begin{document}

\maketitle


\begin{abstract}
       We study the fundamental question of how likely it is that two randomly chosen trees are isomorphic to each other for different models of random trees. We show that the probability decays exponentially for rooted labeled trees as well as for Galton--Watson trees with bounded degrees but that this is not true for plane trees, thus providing a counterexample in the general case of Galton--Watson trees without degree restrictions. We also derive limiting distributions for some related parameters: the number of vertices of given degrees in pairs of labeled trees conditioned on being isomorphic, thus showing that they have a different shape than usual labeled trees, as well as the number of labelings and plane representations of Pólya trees. The results rely on both probabilistic and analytic tools.
\end{abstract}

\section{Introduction}

The study of random trees is a rich and active subject in combinatorial probability theory. Just as it is a classical subject within combinatorics to study the symmetries and isomorphism classes of trees, it is natural to ask about the isomorphism classes of random trees. In this paper we answer the fundamental question of how likely it is that two random trees, drawn at random from some family of such trees, are isomorphic to each other. This question has previously been studied for the family of phylogenetic trees, i.e. full binary trees labeled at the leaves, in \cite{MR2582703}. This is equivalent to the case of labeled binary trees (all full binary trees have the same number of internal and external vertices so every phylogenetic tree of size $n$ corresponds to the same number of fully labeled trees of size $n$). We will study the problem in a more general setting with a mix of probabilistic and analytic tools, while answering other natural questions along the way.

A \textit{Galton--Watson tree} is a rooted tree obtained by a growth process. The root obtains a number of children according to some random variable $X$ that is assumed to take values in the non-negative integers including 0 and some integer larger than 1. Then, the tree grows by letting any offspring in the tree have children according to the same distribution, independently of all other vertices. We will mainly focus on \textit{conditioned} Galton--Watson trees, i.e. trees conditioned on having size $n$. Galton--Watson trees are a special case of (and, under mild conditions, equivalent to) \textit{simply generated trees}, a class of trees that is defined in terms of generating functions. Each simply generated tree has a \textit{weight}, denoted by $w(T)$. Let 
\begin{equation*}
    T(x) = \sum_{T\in\mathcal{T}} w(T) x^{|T|} 
\end{equation*}
be the generating function for a class of such trees and let 
\begin{equation*}
    \Phi(z) = \sum_{k\geq 0} w_k x^k 
\end{equation*}
be the \textit{weight generating function} where $w_k$ is a weight associated with a vertex having $k$ children. Then, $T(x)$ satisfies 
\begin{equation}
\label{eq:SimpleFuncEq}
    T(x) = x\Phi(T(x)) .
\end{equation}
Given the offspring distribution $X$ (or, equivalently, the weight generating function $\Phi(z)$) we let $\mathcal{T}$ denote Galton--Watson (or simply generated) trees and $\mathcal{T}_n$ denote Galton--Watson (or simply generated) trees conditioned on having size $n$. For specific realizations of these types of trees we use $T,T_n$. To get a random tree of size $n$ from a simply generated class we pick one according to its weight
\begin{equation*}
    \Prob(\mathcal{T}_n = T) = \frac{w(T)}{\sum_{S\in\mathcal{T}_n} w(S)} = \frac{w(T)}{[x^n]T(x)}.
\end{equation*}
By choosing different distributions or weights, we can obtain different classes of random trees. Some examples of Galton--Watson and simply generated trees are labeled trees, plane trees and binary trees. The book \cite{DrmotaBook_MR2484382} is a general reference on the topic of random trees. 

The connection between simply generated trees and Galton--Watson trees is that, for the latter, the weight $w_k$ is  taken to be the probability of a vertex having $k$ children. Under the (mild) assumption that there exists a positive $\tau$ within the radius of convergence of $\Phi(z)$ such that
\begin{equation*}
    \Phi(\tau) = \tau\Phi'(\tau) < \infty, 
\end{equation*}
we can always assume that our trees, whether they are conditioned Galton--Watson or simply generated ones, are critical Galton--Watson trees, with slight modifications to the offspring distribution \cite[Subsection 1.2.7]{DrmotaBook_MR2484382}.  

Under the condition that we can find a $\tau$ as above, we can find $\rho=\frac{\tau}{\Phi(\tau)}$ such that \cite[Theorem 3.6]{DrmotaBook_MR2484382}
\begin{equation*}
    [x^n]T(x) \sim C n^{-3/2}\rho^{-n} .
\end{equation*}
Note that for Galton--Watson trees, $[x^n]T(x)$ is exactly the probability  $\Prob(|\mathcal{T}| = n)$. The asymptotic behavior is obtained using the method of singularity analysis, and $\rho$ is the smallest (positive) singularity of the function $T(x)$. It is a general fact that the singularity of a Galton--Watson tree satisfies
\begin{align}
\label{eq:GWSystOfEq}
    \tau = \rho \Phi(\tau) \nonumber,\\
    1 = \rho \Phi'(\tau) ,
\end{align}
which, formulated in words, means that the implicit function theorem fails at the point $(\rho,\tau)$. Classes of Galton--Watson trees for which $\E X = 1$ are called \textit{critical}. In this case, we always have $\rho=1$ so that the probability $P(|\mathcal{T}|=n)$ decays like $n^{-3/2}$.

\textit{Pólya trees} $\mathcal{P}$ are rooted, unordered, unlabeled trees. They are not Galton--Watson trees, even though they share many of the same properties. They can be defined by their generating function
\begin{equation*}
    P(x) = \sum_{P\in\mathcal{P}} x^{|P|} ,
\end{equation*}
that satisfies the functional equation \cite[Section 1.2.5]{DrmotaBook_MR2484382}
\begin{equation*}
    P(x) = x\exp\left(\sum_{j\geq0} \frac{P(x^j)}{j}\right).
\end{equation*}
We can also study Pólya trees $\mathcal{P}_D$ with degrees restricted to lie in some set $D$
\begin{equation*}
    P_D(x) = \sum_{P\in\mathcal{P}_D} x^{|P|}.
\end{equation*}

In this paper we will study rooted labeled trees and general Galton--Watson trees with vertex degrees restricted to lie in some finite set $D$. The isomorphism classes of rooted labeled trees correspond exactly to Pólya trees, and the isomorphism classes of Galton--Watson trees with degrees in the set $D$ to Pólya trees with degrees in the same set. We will talk about Pólya trees and the isomorphism classes of labeled or Galton--Watson trees interchangeably. If we fix the type of Galton--Watson tree, we can also talk about the weight of the isomorphism class associated with a given Pólya tree. We will use $W(T)$, for $T$ a Galton--Watson tree, to denote the total weight of the isomorphism class of $T$. By abuse of notation, we will also use $W(P)$, for $P$ a Pólya tree, to denote the weight of the isomorphism class associated with $P$.

By the orbit-stabilizer theorem \cite[Lemma 6.1]{MR2339282}, the number of isomorphism classes of a tree $T$ under the action of the group $G$ is
\begin{equation*}
    \frac{|G|}{|\aut{T}|} .
\end{equation*}
For rooted labeled trees, $G$ will be the symmetric group $S_n$ acting on the labels so that $|G| = n!$. For Galton--Watson trees, there is always an implicit ordering (we can talk about a first child etc.) and the proper action to study is the one that, for each vertex in the tree, permutes its children. Let $V(T)$ denote the vertex set of the tree $T$, $\deg(T)$ its root degree and $\deg(v)$ the outdegree of a vertex $v\in V(T)$. Also, let $\mathcal{B}(T)$ be the root branches of $T$ and $\mathcal{B}_I(T)$ be the unique root branches up to isomorphism. For the group action on Galton--Watson trees that permutes the branches at every vertex we have
\begin{equation*}
    |G| = \deg(T)! \prod_{B\in\mathcal{B}(T)} |G_{B}| =  \prod_{v\in V(T)} \deg(v)! ,
\end{equation*}
where $G_B$ is the group action restricted to the root branch $B$. In a similar vein, the size of the automorphism group of any rooted tree $T$ satisfies
\begin{equation}
\label{eq:RecAutGroup}
    |\aut{T}| = \prod_{B\in\mathcal{B}_I(T)} \mathrm{mult}(B_i)! |\aut{T_i}|^{\mathrm{mult}(B_i)} ,
\end{equation}
with $\mathrm{mult}(B)$ being the multiplicity of the branch $B$ (which is only taken up to isomorphism). In both of these cases, the expression mirrors the fact that the group is built up by iterated direct and wreath products of symmetric groups (see \cite{MR1373683}, Proposition 1.15).
With the orbit-stabilizer theorem in mind, the total number of plane representations $\mathcal{PR}(T)$ of a Pólya tree $P$ (or an isomorphism class of Galton--Watson trees) is
\begin{equation*}
    \mathcal{PR}(P) = \frac{\prod_{v\in V(P)} \deg(v)!}{|\aut{P}|} .
\end{equation*}
To get the total weight $W(P)$ of the isomorphism class, we need to multiply by the weights:
\begin{equation*}
    W(P) = w(T)\mathcal{PR}(P) = \frac{\prod_{v\in V(P)} w_{\deg(v)}\deg(v)!}{|\aut{P}|} ,
\end{equation*}
where $T$ is some Galton--Watson tree in the isomorphism class corresponding to $P$. 

Let $p_n$ be the probability that two rooted labeled trees on $n$ vertices are isomorphic when we pick them uniformly at random. Then there are
\begin{equation*}
    \sum_{P\in\mathcal{P}_n} \frac{n!^2}{|\aut{P}|^2} ,
\end{equation*}
pairs of isomorphic labeled trees, out of $n^{2(n-1)}$ pairs of such trees in total. We can define the bivariate generating function
\begin{equation*}
    P(x,t) = \sum_{P\in\mathcal{P}} \frac{x^{|P|}}{|\aut{P}|^t}  .
\end{equation*}
Then, the probability is given by
\begin{equation}
\label{eq:probRep}
    p_n = \frac{1}{n^{2(n-1)}} \sum_{P\in\mathcal{P}_n} \frac{n!^2}{|\aut{P}|^2} = \left(\frac{n!}{n^{n-1}}\right)^2 [x^n]P(x,2).
\end{equation}
In other words, we sum over all isomorphism classes of labeled trees (of size $n$) and square the probability of a tree belonging to that class.

If we fix a type of Galton--Watson trees with weights $w_k$ for $k$ in a finite set $D$ (in vector format we can write $\mathbf{w} =[w_1,w_2,\ldots,w_d]$), then we can define the following bivariate generating function for Pólya trees with degrees in $D$
\begin{equation*}
    P_D(x,t) = \sum_{P\in\mathcal{P}_D} W(P)^t x^{|P|} = \sum_{P\in\mathcal{P}_D} \left(\frac{\prod w_kk!}{|\aut{P}|}\right)^t x^{|P|} .
\end{equation*}
Note that the case of $t=1$ gives back $[x^n]T(x)$, the generating function for the corresponding simply generated tree. To study the probability that two trees are isomorphic, we set $t=2$ and divide by $([x^n]T(x))^2$.

\subsection{Results}
For labeled trees we prove that the probability that two trees are isomorphic is, asymptotically, exponentially small. 

\begin{theorem}
\label{thm:problabeled}
The probability $p_n$ that two labeled rooted trees are isomorphic has the full asymptotic expansion
\begin{equation*}
    p_n \sim A n^{3/2} c_l^{n} \left( 1 + \sum_{k=1}^\infty \frac{e_k}{n^k}\right),
\end{equation*}
where $A\approx 2.397678$, $c_l\approx 0.354379$ and the $e_k$ are constants that can be calculated. 
\end{theorem}

For Galton--Watson trees with bounded degrees, we can also prove exponential decay of the probabilities but we have not been able to obtain a full asymptotic expansion except in a special case studied in Subsection \ref{subsec:UnaryBinaryTrees}.

\begin{theorem}
\label{thm:probGW}
The probability $g_n$ that two Galton--Watson trees with degrees in a finite set $D$ are isomorphic satisfies
\begin{equation*}
    g_n \leq B  c_g^{n},
\end{equation*}
for some constants $B$ and $c_g<1$.
\end{theorem}

These results are in line with what has previously been observed for phylogenetic trees \cite{MR2582703} (see also \cite{https://doi.org/10.1002/rsa.20428} for a related question in the context of binary search trees), and also with what one might suspect since this is fundamentally a question of how the weight of simply generated trees, which grows exponentially, relates to their isomorphism classes, which are equivalent to Pólya trees and therefore exponentially many. It could, however, be the case that the weight concentrates on only one or a few isomorphism classes in such a way that the probabilities decay slower than any exponential function. Somewhat surprisingly, we show that this can indeed happen in the general case (i.e. when we do not assume bounded degrees of the Galton--Watson trees). To prove this, we show that the probability that two plane trees are isomorphic exhibits subexponential decay. 

\begin{theorem}
\label{thm:SubexponentialDecay}
The probability that two plane trees are isomorphic decays subexponentially. Thus, we cannot obtain exponential bounds on the probability that two conditioned Galton--Watson trees are isomorphic in general.
\end{theorem}

Since a pair of isomorphic trees will necessarily have the same number of vertices of each degree, we can condition on the event that two labeled trees are isomorphic and study the number of vertices with degree $d$ occurring in one of the trees. The primary motivation is to see if trees conditioned on being isomorphic exhibit a different structure than regular labeled trees. We will see that this is indeed the case when we compare the number of leaves in the two settings in Subsection \ref{subsec:Leaveslabeled}. In the general case, we have the following result on the degree distribution of the vertices.

\begin{theorem}
\label{thm:cltlabeled}
Let $\mathbf{X}_n$ be a random vector that counts the number of vertices of (out)degree $\mathbf{d} = (d_1,d_2,\ldots,d_k)$ in either of a pair of isomorphic labeled trees of size $n$. Then
\begin{align*}
    \E \mathbf{X}_{n} &= \boldsymbol{\mu} n + O(1) , \\
    \cov \mathbf{X}_{n} &= \mathbf{\Sigma} n + O(1) ,
\end{align*}
for a vector $\boldsymbol{\mu}=(\mu_1,\mu_2,\ldots,\mu_k)$ and a matrix $\mathbf{\Sigma} = (\sigma_{i,j})_{1\leq i,j\leq k}$. Furthermore, we have joint convergence to a normal distribution
\begin{equation*}
    \frac{X_{n} - \E\mathbf{X}_n}{\sqrt{n}} \xrightarrow{d} N(\mathbf{0},\mathbf{\Sigma}) .
\end{equation*}
\end{theorem}

If we fix a class of simply generated trees with bounded degrees, the family of Pólya trees with the same degree restrictions corresponds to the isomorphism classes of simply generated trees. It is then possible to ask what the weight of a randomly chosen Pólya tree is. This is particularly interesting when the weights $w_i$ are integers and the coefficients of the generating function correspond to the number of trees. This is for example the case for binary, unary-binary and ternary trees. In these cases, the weight of a given isomorphism class is equal to the number of plane representations of the corresponding unordered, unlabeled tree. Thus, the following theorem shows that the logarithm of the number of plane representations of a randomly chosen Pólya tree with degree restrictions satisfies a central limit theorem. 

\begin{theorem}
\label{thm:nrRepGW}
    Let $\mathcal{P}_n$ be a random tree of size $n$ from the class of unordered, unlabeled trees with degrees in the finite set $D$ and $\mathcal{T}_D$ be a class of Galton--Watson trees with the same degree restrictions. Then the weight $W(\mathcal{P}_n)$ of $\mathcal{P}_n$ seen as an isomorphism class of $\mathcal{T}_D$ has expected value and variance
    \begin{align*}
        \E [\log W(\mathcal{P}_n)] &= \mu n + O(1) , \\
        \Var [\log W(\mathcal{P}_n)] &= \sigma^2 n + O(1) ,
    \end{align*}
    for some constants $\mu$ and $\sigma$ and satisfies
    \begin{equation*}
        \frac{\log W(\mathcal{P}_n) - \mu n}{\sqrt{n}} \xrightarrow[]{d} N(0,\sigma^2) .
    \end{equation*}
\end{theorem}

We can also derive an analogue for the set of labeled trees without degree restrictions. 
\begin{theorem}
\label{thm:nrRepLab}
    Let $\mathcal{P}_n$ be a random Pólya tree of size $n$, then the number of labelings $\operatorname{L}(\mathcal{P}_n)$ of $\mathcal{P}_n$ has expected value and variance
    \begin{align*}
        \E [\log\operatorname{L}(\mathcal{P}_n)] &= n\log n - (\mu+1) n + \frac{\log n}{2}+ O(1) , \\ 
        \Var [\log\operatorname{L}(\mathcal{P}_n)] &= \sigma^2 n + O(1) ,
    \end{align*}
    for numerical constants $\mu \approx 0.137342$ and $\sigma^2 \approx 0.196770$ and satisfies
    \begin{equation*}
        \frac{\log \operatorname{L}(\mathcal{P}_n) - \E [\log\operatorname{L}(\mathcal{P}_n)]}{\sqrt{n}} \xrightarrow[]{d} N(0,\sigma^2) .
    \end{equation*}
\end{theorem}
Note that the mean is of the order $n\log n$ which is somewhat unusual for combinatorial limit laws. It is worth pointing out that we have similar results for labeled trees with degree restrictions by Theorem \ref{thm:nrRepGW} above, but then both the mean and variance are of order $n$.

We make a passing remark that both the number of vertices of given degrees and the logarithm of the weight of isomorphism classes are examples of so called \textit{additive functionals}. An additive functional is a function $F(T)$ of rooted trees that satisfies the recursion
\begin{equation*}
    F(T) = f(T) + \sum_{i=1}^r F(S_i) , 
\end{equation*}
where we sum over the $r$ root branches $S_1,S_2,\ldots,S_r$ and $f(T)$ is called the toll function associated with the functional. Properties of additive functionals have been studied extensively for various types of random trees, see \cite{MR3318048, MR3432572, almostLocal_2020} (and further references therein) for some examples involving Galton--Watson and Pólya trees.

In Section \ref{sec:funcEq} we derive functional equations to be used in the other sections. We then prove Theorem \ref{thm:problabeled} in Section \ref{sec:labeled} and Theorem \ref{thm:probGW} in Section \ref{sec:GWBounded} while we provide a counterexample that proves Theorem \ref{thm:SubexponentialDecay} in Section \ref{sec:counterex}. Section \ref{sec:clt} is devoted to proving Theorem \ref{thm:cltlabeled} and Section \ref{subsec:nrRep} to proving Theorems \ref{thm:nrRepGW} and \ref{thm:nrRepLab}. We also provide examples, estimating the probability that two  unary-binary trees are isomorphic in Subsection \ref{subsec:UnaryBinaryTrees}, estimating the mean number of leaves in a pair of isomorphic labeled trees in Subsection \ref{subsec:Leaveslabeled} and finding numerical values for the moments in the central limit theorem for the number of plane representations of Pólya trees where all vertices have 0, 1 or 2 vertices in Subsection \ref{subsec:EstPlaneRepBin}.

\section{Functional equations}
\label{sec:funcEq}

The functional equation for the bivariate generating function counting the size of the automorphism group of Pólya trees has previously been derived in \cite{autRandTrees}. We reproduce the calculations here for completeness. 

We have the following symbolic decomposition of Pólya trees
\begin{equation*}
    \mathcal{P} = \bullet \times \bigotimes_{P\in\mathcal{P}} (\emptyset \uplus \{P\} \uplus \{P,P\} \uplus \ldots )  ,
\end{equation*}
which, by (\ref{eq:RecAutGroup}), translates to the functional equation
\begin{equation}
\label{eq:PolyaFuncEqProd}
    P(x,t) =  x \prod_{P\in\mathcal{P}}\left(\sum_{n=0}^\infty x^{n|P|} n!^{-t} |\aut{P}|^{-nt}\right)   
\end{equation}
when we take automorphisms into account. We rewrite this as follows:
\begin{align}
\label{eq:PolyaFuncEq1}
&P(x,t) =  x \exp\left(\sum_{P\in\mathcal{P}} \log \left(\sum_{n=0}^\infty x^{n|P|} n!^{-t} |\aut{P}|^{-nt}\right)\right) \nonumber \\
&= x \exp\left(\sum_{P\in\mathcal{P}} \sum_{k=1}^\infty \frac{(-1)^{k-1}}{k} \left(\sum_{n=1}^\infty x^{n|P|} n!^{-t} |\aut{P}|^{-nt}\right)^k \right) \nonumber \\
&= x \exp\Big(\sum_{P\in\mathcal{P}} \sum_{k=1}^\infty \frac{(-1)^{k-1}}{k} \sum_{\substack{\lambda_1+\lambda_2\\+\ldots = k}} \binom{k}{\lambda_1,\lambda_2,\ldots} \prod_{n=1}^\infty \big(x^{n|P|} n!^{-t} |\aut{P}|^{-nt}\big)^{\lambda_n} \Big).
\end{align}
For an integer partition $\lambda$, write $\lambda = (\lambda_1,\lambda_2,\ldots)$, where $\lambda_i$ is the number of $i$'s in the partition. We let $|\lambda| = \lambda_1+\lambda_2+\ldots$ denote the total number of summands, and we write $\lambda \vdash j$ to denote that $\lambda$ is a partition of $j$, i.e., $j = \lambda_1 + 2\lambda_2+3\lambda_3 + \ldots$. We use this to rewrite the expression inside the exponential function in (\ref{eq:PolyaFuncEq1}).
\begin{multline*}
    \sum_{P\in\mathcal{P}} \sum_{k=1}^\infty \frac{(-1)^{k-1}}{k} \sum_{j=1}^\infty  
    \sum_{\substack{\lambda_1+\lambda_2+\ldots = k\\\lambda_1+2\lambda_2+\ldots = j}} \binom{k}{\lambda_1,\lambda_2,\ldots}  x^{j|P|}  |\aut{P}|^{-jt}\prod_{n=1}^\infty n!^{-\lambda_n t} \\
    = 
    \sum_{j=1}^\infty \sum_{\lambda\vdash j}  \frac{(-1)^{|\lambda|-1}}{|\lambda|}  \binom{|\lambda|}{\lambda_1,\lambda_2,\ldots} \left(\prod_{n=1}^\infty n!^{-\lambda_n t}\right) \sum_{P\in\mathcal{P}} x^{j|P|}  |\aut{P}|^{-jt}  \\
    = 
    \sum_{j=1}^\infty \sum_{\lambda\vdash j}  \frac{(-1)^{|\lambda|-1}}{|\lambda|}  \binom{|\lambda|}{\lambda_1,\lambda_2,\ldots} \left(\prod_{n=1}^\infty n!^{-\lambda_n t}\right)  P(x^j,jt) .
\end{multline*}
We now have the functional equation
\begin{equation*}
    P(x,t) = x\exp\left( P(x,t) + \sum_{j=2}^\infty \frac{c(j,t)}{j} P(x^j,jt)\right) ,
\end{equation*}
where we define
\begin{equation*}
    c(j,t) = j\sum_{\lambda \vdash j} \frac{(-1)^{|\lambda|-1}}{|\lambda|} \binom{|\lambda|}{\lambda_1, \lambda_2,\ldots }\left( \prod_{m=1}^\infty m!^{-\lambda_m t} \right) .
\end{equation*}
We recover the case that we are interested in by setting $t=2$, which yields
\begin{equation}
\label{eq:PolyaGenFun}
    P(x,2) = x\exp\left( P(x,2) + \sum_{j=2}^\infty \frac{c(j,2)}{j} P(x^j,2j)\right) .
\end{equation}
For later reference, we note that we can expand (\ref{eq:PolyaGenFun}) and collect according to root degree as
\begin{equation}
    \label{eq:labeled_expand}
    P(x,2) = x\sum_{k\geq 0} \sum_{\lambda \vdash k} \prod_j \frac{\left(c(j,2)P(x^j,2j)\right)^{\lambda_j}}{j^ {\lambda_j}\lambda_j!} .
\end{equation}
This can be seen by introducing an extra variable $y$ in (\ref{eq:PolyaFuncEqProd}) that keeps track of the root degree as follows
\begin{equation*}
    P(x,t) =  x \prod_{P\in\mathcal{P}}\left(\sum_{n=0}^\infty x^{n|P|} n!^{-t} |\aut{P}|^{-nt} y^n\right)   ,
\end{equation*}
performing the same calculations as above, and then extracting the coefficient in front of $y^k$ to get the $k$-th term.

To obtain a functional equation for the weights of isomorphism classes of Galton--Watson trees with restricted degrees we can perform similar calculations with $W(P)$ instead of $|\aut{P}|^{-1}$. We extract only the terms corresponding to root degrees lying in $D$, and then multiply by the extra factors $(w_k k!)^2$, i.e. the contribution to the weight coming from the root. This gives
\begin{equation*}
    P_D(x,2) = x\sum_{k\in D} (w_k k!)^2 \sum_{\lambda \vdash k} \prod_j \frac{\left(c(j,2)P_D(x^j,2j)\right)^{\lambda_j}}{j^{\lambda_j}\lambda_j!}, 
\end{equation*}
which is a polynomial in the factors $P_D(x^j,2j)$ since $D$ is a finite set. More generally, we have
\begin{equation}
\label{eq:funcEqBoundDeg}
    P_D(x,t) = x\sum_{k\in D} (w_k k!)^t \sum_{\lambda \vdash k} \prod_j \frac{\left(c(j,t)P_D(x^j,jt)\right)^{\lambda_j}}{j^{\lambda_j}\lambda_j!} ,
\end{equation}
by the same reasoning.

\section{Singularity analysis for labeled trees}
\label{sec:labeled}

We can now solve the functional equation (\ref{eq:PolyaGenFun}) in terms of the tree function $Y(x)$, which is the exponential generating function of the class of rooted labeled trees and satisfies $Y(x) = xe^{Y(x)}$. We then get
\begin{equation}
\label{eq:labeledSolution}
    P(x,2) = Y\left(x\exp\left(\sum_{j=2}^\infty \frac{c(j,2)}{j} P(x^j,2j)\right)\right) .
\end{equation} 
We want to bound the radius of convergence $\alpha$ of $P(x,2)$. Let $\rho_P \approx 0.338$ be the known singularity of the univariate generating function for Pólya trees $P(x) = P(x,0)$. See \cite[Chapter 3.1]{DrmotaBook_MR2484382} for more background on $Y(x)$ and $P(x)$. Now note that $\alpha$ is at least as large as $\rho_P$ since $|\aut{P}|^{-2}\leq 1$ for all $P$. We can get an upper bound on $\alpha$ by first noting that for the probability $p_n$ of two labeled trees being isomorphic we have, by Cauchy--Schwarz, $p_n \geq  \frac{1}{|\mathcal{P}_n|}$, i.e. we cannot do worse than when all of the isomorphism classes are equally likely. By $(\ref{eq:probRep})$ and Stirling's approximation, the coefficients of $P(x,2)$ are $\sim\frac{e^{2n}}{2\pi n^3}$ times as large 
as $p_n$:

\begin{equation*} 
    [x^n]P(x,2) \sim \frac{e^{2n}}{2\pi n^3} p_n \geq \frac{e^{2n}}{2\pi n^3 |\mathcal{P}_n|}.
\end{equation*}
We know that the radius of convergence of $P(x)$ is $\rho_P$, implying that the radius of convergence of $\sum \frac{1}{\mathcal{P}_n} x^n$ is $1/ \rho_P$. In other words, $\frac{1}{\mathcal{P}_n}$ decays exponentially with base $\rho_P$. Since larger coefficients imply a smaller radius of convergence, and since only the exponential factors affect the radius, this gives an upper bound of
\begin{equation*}
\alpha \leq 1/(\rho_P e^2) \approx 0.400 .
\end{equation*}
As this is smaller than 1, we find that the radius of convergence of $P(x^j,2)$ is larger than that of $P(x,2)$ for any $j\geq 2$. Furthermore, 
\begin{equation*}
    \frac{1}{|\aut{P}|^{2j}} \leq \frac{1}{|\aut{P}|^2} 
\end{equation*}
which means that the coefficients of $P(x^j,2j)$ are smaller than those of $P(x^j,2)$. As a consequence, $P(x^j,2j)$ is analytic in a larger region than $P(x^j,2)$ and, by extension, $P(x,2)$.

From the definition of $c(j,t)$, we get that
\begin{equation*}
     \frac{|c(j,2)|}{j} \leq \sum_{\lambda \vdash j} \binom{|\lambda|}{\lambda_1, \lambda_2,\ldots } .
\end{equation*}
This counts the number of partitions of $j$ when we take order of the parts into account. By definition this is the number of compositions of $j$, which is known to be $2^{j-1}$. Thus, we have
\begin{multline*}
    \left|\sum_{j=2}^\infty \frac{c(j,2)}{j} P(x^j,2j)\right| \leq \sum_{j=2}^\infty 2^j \sum_P \frac{|x|^{j|P|}}{|\aut{P}|^{2j}} \\
    \leq \sum_P \frac{1}{|\aut{P}|^{2}} \sum_{j=2}^\infty 2^j |x|^{j|P|} = \sum_P \frac{4|x|^{2|P|}}{|\aut{P}|^{2}(1-2|x|^{|P|})} = O(P(|x|^{2},2)) .
\end{multline*}
Now, the Weierstrass M-test gives that
\begin{equation*}
    \sum_{j=2}^\infty \frac{c(j,2)}{j} P(x^j,2j) 
\end{equation*}
is a uniformly convergent sum of analytic functions in some region $|x|<\alpha+\epsilon$ for $\epsilon>0$, so that the sum, too, is analytic there. If we let $\xi(x)$ denote the expression inside the tree function in (\ref{eq:labeledSolution}), then it is the composition and product of analytic functions so that it, too, is an analytic function.

We can compose the singular expansion of $Y(x)$ at its dominant singularity with the Taylor expansion of $\xi(x)$ at $\alpha$. Since $Y(x)$ is the exponential generating function of labeled trees, which are simply generated, its singular expansion can be found to be
\begin{equation*}
    Y(x) = 1 - \sqrt{2}\sqrt{1-ex} + \ldots ,
\end{equation*}
by general principles \cite{DrmotaBook_MR2484382}, while the dominant singularity can be seen to be $\frac{1}{e}$. We find that
\begin{equation*}
    P(x,2) = 1 - \sqrt{2}\sqrt{1-e(\xi(\alpha) + \xi'(\alpha)(x-\alpha) + \frac{1}{2}\xi''(\alpha)(x-\alpha)^2 + \ldots)} + \ldots ,
\end{equation*}
and we see that $\alpha$ must be the number that satisfies $\xi(\alpha) = e^{-1}$ for it to be a singularity. Using this fact together with the Taylor expansion for $\sqrt{1+z}$ means that the expression is equal to
\begin{equation*}
    1 - \sqrt{2\alpha e\xi'(\alpha)}\sqrt{1-\frac{x}{\alpha}}\left(1 - \frac{\alpha \xi''(\alpha)}{2\xi'(\alpha)}  \left(1-\frac{x}{\alpha}\right) + \ldots \right) .
\end{equation*}
A standard application of singularity analysis shows that, asymptotically, the coefficients behave like 
\begin{equation*}
    [x^n]P(x,2) \sim \sqrt{\frac{e\alpha\xi'(\alpha)}{2\pi}} \frac{\alpha^{-n}}{n^{3/2}} + \ldots, 
\end{equation*}
and by plugging this into \eqref{eq:probRep} and using Stirling's approximation for $n!$, we get
\begin{equation*}
    p_n \sim \sqrt{2\pi e \alpha \xi'(\alpha)} n^{3/2}\frac{1}{(e^2\alpha)^n} + \ldots,
\end{equation*}
for the asymptotic probability that two trees are isomorphic. If we set $A = \sqrt{2\pi e \alpha \xi'(\alpha)}$ and $c_l=\frac{1}{(e^2\alpha)}$, we obtain Theorem \ref{thm:problabeled}.

We can obtain a numerical estimate of $\alpha$ by truncating the series
\begin{equation*}
    \sum_{j=2}^\infty \frac{c(j,2)}{j} P(x^j,2j)
\end{equation*}
as well as the power series $P(x^j,2j)$ in each of its terms. Based on this, we can estimate the smallest positive root of $\xi(x)-e^{-1}$. Then, an estimate of $\xi'(\alpha)$ is obtained by truncating the series and plugging in the estimate for $\alpha$. This gives $A  \approx 2.397678$ and $c_l \approx 0.354379$. 


\section{Galton--Watson trees with bounded degrees}
\label{sec:GWBounded}
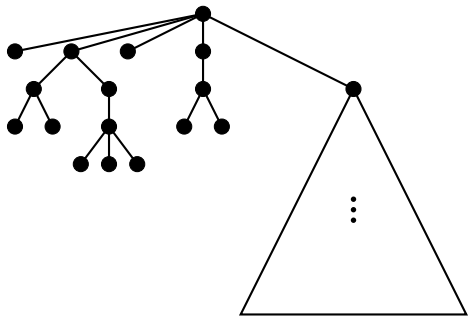
\begin{figure}
    \centering
    
    \begin{tikzpicture}
    {

    \draw[fill] (1,5) circle (0.1);

    \foreach \x in {-1.5,-0.75,0,1}
    {
        \draw[fill] (\x,4.5) circle (0.1);
        \draw[thick] (\x,4.5) -- (1,5);
    }

    \draw[fill] (-0.75-0.5,4) circle (0.1);
    \draw[thick] (-0.75-0.5,4) -- (-0.75,4.5);
    \draw[fill] (-0.75+0.5,4) circle (0.1);
    \draw[thick] (-0.75+0.5,4) -- (-0.75,4.5);

    \draw[fill] (1,4) circle (0.1);
    \draw[thick] (1,4) -- (1,4.5);

    \draw[fill] (-0.75-0.75,3.5) circle (0.1);
    \draw[thick] (-0.75-0.75,3.5) -- (-0.75-0.5,4);
    \draw[fill] (-0.75-0.25,3.5) circle (0.1);
    \draw[thick] (-0.75-0.25,3.5) -- (-0.75-0.5,4);

    \draw[fill] (-0.75+0.5,3.5) circle (0.1);
    \draw[thick] (-0.75+0.5,3.5) -- (-0.75+0.5,4);

    \draw[fill] (0.75,3.5) circle (0.1);
    \draw[thick] (0.75,3.5) -- (1,4);
    \draw[fill] (1.25,3.5) circle (0.1);
    \draw[thick] (1.25,3.5) -- (1,4);

    \draw[fill] (-0.75+0.125,3) circle (0.1);
    \draw[thick] (-0.75+0.125,3) -- (-0.75+0.5,3.5);
    \draw[fill] (-0.75+0.5,3) circle (0.1);
    \draw[thick] (-0.75+0.5,3) -- (-0.75+0.5,3.5);
    \draw[fill] (-0.75+0.875,3) circle (0.1);
    \draw[thick] (-0.75+0.875,3) -- (-0.75+0.5,3.5);

    \draw[fill] (3,4) circle (0.1);
    \draw[thick] (3,4) -- (1,5);
    \draw[thick] (3+1.5,1)-- (3,4) -- (3-1.5,1) -- cycle;
    \node at (3,2.5) {\Huge $\vdots$};
    
    }
    \end{tikzpicture}\\
    \caption{A visual example of a Galton--Watson tree with a giant branch.}
    \label{fig:giant}
\end{figure}

We now show that we have exponential decay for the probability that two Galton--Watson trees with bounded degrees are isomorphic. For an isomorphism class $I$, say that a conditioned Galton--Watson tree has property $\mathcal{S}_I$ if it consists of a tree from $I$ together with a giant branch, i.e. a branch that contains the rest of the $n-|I|$ vertices, attached to the root of $I$ (see Figure \ref{fig:giant}). 
We know that a conditioned Galton--Watson tree has such a giant branch with probability tending to 1 \cite[Section 7]{JansonAllocConden}. Furthermore $P(\mathcal{T}_n\textrm{ has property }\mathcal{S}_I)$ converges to some probability $p_I$ such that
\begin{equation*}
    \sum_I p_I = 1,
\end{equation*}
when we sum over all possible isomorphism classes. To see this, consider a tree of size $n$ with one giant branch of size $n-|T_I|$, where $T_I$ is a representative of the isomorphism class $I$ of the rest of the tree. The limiting probability of $\{\mathcal{T}_n\textrm{ has property }\mathcal{S}_I\}$ is 
\begin{equation*}
    \frac{\Prob(|\mathcal{T}|=n-|T_I|)}{\Prob(|\mathcal{T}|=n)} \mathcal{PR}(I) (\deg(T_I)+1) w_{\deg(T_I)+1} \prod_{B\in \mathcal{B}(T_I)} w(B) ,
\end{equation*}
as we need to choose a tree of size $n-|T_I|$ for the giant branch, an embedding of $I$, as well as a position to attach the giant branch to $I$. Recall that we have $\Prob(|\mathcal{T}|=n)\sim C n^{-3/2}\rho^{-n}$, so, in the limit, this probability is 
\begin{equation*}
    \rho^{|T_I|} \mathcal{PR}(I) (\deg(T_I)+1) w_{\deg(T_I)+1} \prod_{B\in \mathcal{B}(T_I)} w(B) .
\end{equation*}
We can perform calculations where we let $n\to\infty$ and sum over all $I$ to obtain
\begin{multline*}
    \sum_I \rho^{|T_I|} \mathcal{PR}(I) (\deg(T_I)+1) w_{\deg(T_I)+1} \prod_{B\in \mathcal{B}(T_I)} w(B) \\
    = \sum_{S\in\mathcal{T}} \rho^{|S|} (\deg(S)+1) w_{\deg(S)+1} \prod_{B\in \mathcal{B}(S)} w(B) .
\end{multline*}
If we rewrite the expression so that we instead sum over the root degree of $S$, we obtain
\begin{multline*}
     \sum_{k\geq 0} (k+1) w_{k+1} \rho\left( \sum_{B_1\in\mathcal{T}} \rho^{|B_1|} w(B_1)\right)\cdot\ldots\cdot\left( \sum_{B_k\in\mathcal{T}} \rho^{|B_k|} w(B_k)\right)\\ 
    = \rho \sum_{k\geq0} (k+1) w_{k+1} T(\rho)^k = \rho \Phi'(T(\rho)) = 1,
\end{multline*} 
where we use the general fact that $\rho\Phi'(\tau) = 1$ (recall equation \eqref{eq:GWSystOfEq}) in the last step. As the sum is 1, we indeed have a giant branch and the associated probability distribution when $n\to\infty$.

Let $d$ be the maximal degree that a vertex can have. Fix $0<\epsilon<1$ and choose a finite family of isomorphism classes $\mathcal{F}$ as well as a positive integer $n_0$ such that, for all $n\geq n_0$, 
\begin{multline}
\label{eq:ExpBoundsAssumption}
    \sum_{I\in\mathcal{F}} \Prob(\mathcal{T}_n\textrm{ has property }\mathcal{S}_I)^2 \\
    + d!\Prob(\mathcal{T}_n\textrm{ does not have property }\mathcal{S}_I\textrm{ for any }I) < 1-\epsilon .
\end{multline}
This can be done since
\begin{equation*}
    \sum_{I\in\mathcal{F}} p_I^2 < \left(\sum_{I\in\mathcal{F}} p_I\right)^2 \leq 1 ,
\end{equation*}
since we can safely assume that $p_I>0$ for more than one $I$. At the same time, we can take the family $\mathcal{F}$ to be large enough to make 
\begin{equation*}
    d! \sum_{I\notin\mathcal{F}} p_I 
\end{equation*}
arbitrarily small. As the probability $\Prob(\mathcal{T}_n\textrm{ has property }\mathcal{S}_I)$ converges to $p_I$, and $\Prob(\mathcal{T}_n\textrm{ does not have property }\mathcal{S}_I\textrm{ for any }I)$ behaves like 
\begin{equation*}
    \sum_{I\notin\mathcal{F}} p_I 
\end{equation*}
in the limit, the inequality (\ref{eq:ExpBoundsAssumption}) must hold for large enough $n$.

Take some constant $0<c<1$, to be determined later. We want to show that $\Prob(\mathcal{T}_n^{(1)} \simeq \mathcal{T}_n^{(2)})\leq c^{n-n_0}$ for all $n$. As the basis for our induction argument, it is clear that $\Prob(\mathcal{T}_n^{(1)} \simeq \mathcal{T}_n^{(2)})\leq c^{n-n_0}$ for all $n\leq n_0$ as the bound is trivially larger than 1 in that case. Now assume that the induction hypothesis holds and take two conditioned Galton--Watson tree $\mathcal{T}_n^{(1)}$ and $\mathcal{T}_n^{(2)}$ of size $n$. 
If $\mathcal{T}_n^{(1)}$ has property $\mathcal{S}_I$, then the two trees can only be isomorphic if $\mathcal{T}_n^{(2)}$ does as well. If they both have property $\mathcal{S}_I$ for some $I\in\mathcal{F}$, we use the induction hypothesis on the giant branch, and otherwise we use the induction hypothesis on each of the, at most $d$, root branches. If, in the latter case, there are $k$ root branches having sizes $m_1,m_2,\ldots,m_k$, then we get 
\begin{equation*}
    c^{m_1-n_0 + m_2-n_0 + \ldots + m_k-n_0} = c^{n-1-kn_0} \leq c^{n-1-dn_0}
\end{equation*}
as an upper bound on the probability that the two, original, trees are isomorphic. Summing over all alternatives, we obtain the bound
\begin{multline*}
    \Prob(\mathcal{T}_n^{(1)} \simeq \mathcal{T}_n^{(2)}) 
    \leq \sum_{I\in\mathcal{F}} \Prob(\mathcal{T}_n\textrm{ has property }\mathcal{S}_I)^2 c^{n-|I|-n_0} \\
    + \Prob(\mathcal{T}_n\textrm{ does not have property }\mathcal{S}_I\textrm{ for any }I) d! c^{n- dn_0 - 1} .
\end{multline*}
where we get a factor of $d!$ since any permutation of the root branches of a tree $\mathcal{T}_n$ yields a tree isomorphic to it.
We are free to assume that $|I|$ is bounded by $n_0$ and, therefore, also by $1+(d-1)n_0$. This leads to a bound on the expression above in the form of
\begin{multline*}
    c^{n-n_0} \bigg(\sum_{I\in\mathcal{F}} \Prob(\mathcal{T}_n\textrm{ has property }\mathcal{S}_I)^2  \\
    + d!\Prob(\mathcal{T}_n\textrm{ does not have property }\mathcal{S}_I\textrm{ for any }I \in \mathcal{F})  \bigg) c^{-(1 + (d-1)n_0)}\\
    \leq c^{n-n_0}(1-\epsilon)c^{-(1 + (d-1)n_0)} .
\end{multline*}
Note that we can take $c$ arbitrarily close to 1 to guarantee that 
\begin{equation*}
(1-\epsilon)c^{-(1 + (d-1)n_0)}\leq 1, 
\end{equation*}
so that we have the bound
\begin{equation*}
    \Prob(\mathcal{T}_n^{(1)} \simeq \mathcal{T}_n^{(2)}) \leq c^{n-n_0} .
\end{equation*}
This finishes the inductive step and proves Theorem \ref{thm:probGW} with $B=c^{-n_0}$ and $c_g=c$. The argument indicates that the correct value of asymptotic decay might be the constant $c$ that satisfies
\begin{equation*}
    \sum_{I\in\mathcal{P}} p_I^2 c^{-|I|} = 1 ,
\end{equation*}
provided such a constant exists.

\subsection{Singularity analysis for unary-binary trees}
\label{subsec:UnaryBinaryTrees}

One could hope to derive asymptotics for the probabilities through the functional equations obtained in Section \ref{sec:funcEq} and singularity analysis, but the necessary calculations turn out to be technical and we have not found a way to make the argument rigorous. We can, however, study the special case of unary-binary trees, having the weights $w_0=1,w_1=1$, and $w_2=1$, using a specialized method. Note also that the case of phylogenetic trees was studied in \cite{MR2582703}, and that this extends to exponential decay of probabilities for binary trees with weights $\mathbf{w}=[1,0,1]$ and, by extension through a well known bijection, to trees with weights $\mathbf{w}=[1,2,1]$ as well.

For this type of trees, we obtain, from (\ref{eq:funcEqBoundDeg}), the functional equation 
\begin{equation}
\label{eq:UBtreest}
    u(x,t) = x + xu(x,t) + 2^{t-1} x u(x,t)^2 + (1- 2^{t-1}) u(x^2,2t) ,
\end{equation}
where we let $u(x,t) = P_{\{0,1,2\}}(x,t)$ for convenience. We are interested in the case when $t=2$, for which we have
\begin{equation}
\label{eq:UBtreest2}
    u(x,2) = x + xu(x,2) + 2 x u(x,2)^2 - u(x^2,4) .
\end{equation}
Note, however, that when $t=1$ we recover the plane unary-binary trees and, as this is a well studied class of tree, we know that the radius of convergence of $u(x,1)$ is $1/3$. This can also be seen by solving equation \eqref{eq:UBtreest} with $t=1$ to get
\begin{equation*}
    u(x,1) = \frac{1}{2x}\left(1-x-\sqrt{1-2x-3x^2}\right) 
\end{equation*}
and finding the positive $x$ that makes the expression under the square root equal to zero. Since $u(x,2)$ dominates $u(x,1)$ term-wise, its radius of convergence must be at most 1/3.

On the other hand, setting $t=4$ yields
\begin{equation*}
    u(x,4) = x + xu(x,4) + 8 x u(x,4)^2 - 7 u(x^2,8) ,
\end{equation*}
and by recursively calculating the coefficients we see that $u(x,4)$ is term-wise dominated by the solution to
\begin{equation*}
    v(x) = x + xv(x) + 8 x v(x)^2 .
\end{equation*}
The solution to this functional equation that has positive coefficients is 
\begin{equation*}
    v(x) = \frac{1-x-\sqrt{1-2x-31x^2}}{16x} ,
\end{equation*}
and its radius of convergence, readily calculated to be $\frac{4\sqrt{2}-1}{31}$, is a lower bound for the radius of convergence of $u(x,4)$. Since $\sqrt{\frac{4\sqrt{2}-1}{31}}>1/3$, we see that $u(x^2,4)$ is convergent in a larger region than $u(x,2)$. Now, we can apply the method of singularity analysis to find the asymptotic growth of the coefficients. We apply the following general theorem:

\begin{theorem}[Theorem 2.19 from \cite{DrmotaBook_MR2484382}]
\label{thm:singAnalysis}
Suppose that $F(x, y)$ is an analytic function in $x, y$ around
$x = y = 0$ such that $F(0, y) = 0$ and that all Taylor coefficients of $F$ around $0$ are real and non-negative. Then there exists a unique analytic solution $y = y(x)$ of the functional equation
\begin{equation*}
    y=F(x,y) ,
\end{equation*}
with $y(0) = 0$ that has non-negative Taylor coefficients around 0. 

If the region of convergence of $F(x, y)$ is large enough such that there exist
positive solutions $x = x_0$ and $y = y_0$ of the system of equations
\begin{align}
\label{eq:sysEqSingAn}
    y = F(x,y) , \nonumber \\
    1 = F_y(x,y) ,
\end{align}
with $F_x(x_0,y_0)\neq0$ and $F_{yy}(x_0,y_0)\neq0$, then $y(x)$ is analytic for $|x| < x_0$ and there exist functions $g(x)$, $h(x)$ that are analytic around $x = x_0$ such that $y(x)$ has a representation of the form
\begin{equation}
\label{eq:funcEqThm1}
    y(x) = g(x) - h(x)\sqrt{1-\frac{x}{x_0}} ,
\end{equation}
locally around $x=x_0$. We have $g(x_0) = y(x_0)$ and
\begin{equation*}
    h(x_0) = \sqrt{\frac{2x_0F_x(x_0,y_0)}{F_{yy}(x_0,y_0)}}
\end{equation*}
Moreover, (\ref{eq:funcEqThm1}) provides a local analytic continuation of $y(x)$ (for $\arg(x - x_0) \neq  0$).

If we assume that $[x^n] y(x) > 0$ for $n \geq n_0$, then $x = x_0$ is the only singularity of $y(x)$ on the circle $|x| = x_0$ and we obtain an asymptotic expansion for $[x^n] y(x)$ of the form
\begin{equation}
\label{eq:genAsymEsp}
    [x^n]y(x) = \sqrt{\frac{x_0F_x(x_0,y_0)}{2\pi F_{yy}(x_0,y_0)}} x_0^{-n} n^{-3/2} \left(1+O(n^{-1})\right) .
\end{equation}
\end{theorem}

In our case, we see from \eqref{eq:UBtreest2} that
\begin{equation*}
    F(x,y) = x + xy + 2 x y^2 - u(x^2,4) .
\end{equation*}
Checking the conditions, we see that there must exist a point $(x_0,y_0)$ of the kind specified in the theorem since, by the implicit function theorem, the only other options are that the functional equation is itself singular (which we have shown to be false) or that $P_D(x,2)$ has no singularity (which we also know to be false). To obtain the probabilities we seek, we now only need to divide by the square of the asymptotic number of unary-binary trees. 

To estimate the constants in \eqref{eq:genAsymEsp} numerically, we approximate $F(x,y)$ by a simpler expression. This can be done by truncating the, rapidly converging, series $u(x^2,4)$ to get a polynomial $U(x)$. We can then use standard numerical solvers for the system of equations
\begin{equation*}
    \begin{cases}
        y = x + xy + 2 x y^2 - U(x) ,\\
        1 = x + 4xy ,
    \end{cases}
\end{equation*}
which is an approximation of the system \eqref{eq:sysEqSingAn}. This gives us estimates of $x_0,y_0$ and we can, furthermore, estimate $F_x(x_0,y_0)$ and $F_{yy}(x_0,y_0)$ by truncating these expressions as well. Finding (the square) of the number of unary-binary trees can be done by, similarly, solving the system of equations \eqref{eq:sysEqSingAn} based on \eqref{eq:UBtreest} with $t=1$. We find that the probability of two unary-binary trees being isomorphic behaves like
\begin{equation*}
   g_n \sim C n^{3/2}\delta^n 
\end{equation*}
for large $n$, where $C \approx 1.279101 $ and $\delta \approx 0.412681$.


Note that this method does not appear to generalize well, not even to binary trees (with weighs $w_0=1,w_1=2,w_2=1$) despite the fact that these trees have been studied before (though with different methods). The problem is that we get an overlap between the lower and upper bounds that we want to compare. If we instead consider trees with larger out-degrees, using numerical methods, we again appear to have too weak bounds except in somewhat artificial examples such as ternary trees with weights $\textbf{w}=[1,1,1,0.5]$.

\section{Counterexample to exponential decay for general Galton--Watson trees}
\label{sec:counterex}

We prove Theorem \ref{thm:SubexponentialDecay} by creating a counterexample for plane trees, one type of simply generated tree which has the weights $w_k=1$. We use the pigeonhole principle to find a set of different isomorphism classes that all contain many plane trees. By attaching all the isomorphism classes in the set to a common root, we obtain an isomorphism class that contains a large fraction of plane trees of a given size. 

It is well known that the number of plane trees of size $n$ is the $n$-th Catalan number $C_n$ and, thus, asymptotically of order $4^{n-O(\log{n})}$. This implies that the number of plane representations of any one Pólya tree (which we recall to be an isomorphism class of plane rooted trees) is less than $4^n$. We can now divide the interval $[0,\log{(4)} n]$ into $cn$ subintervals for some small $0<c<1$ (this is not guaranteed to be an integer but any rounding error is irrelevant in the limit) and partition $\mathcal{P}_n$ according to which interval $\log{(\mathcal{PR}(P))}$ falls in. Summing all plane representations of Pólya trees in $\mathcal{P}_n$ gives us the number of plane trees of size $n$ so, by the pigeonhole principle, one of the intervals must contain a set of Pólya trees that has 
\begin{equation*}
    \frac{C_n}{cn} = 4^{n-O(\log{n})} 
\end{equation*}
plane representations in total. Furthermore, there is an $x$ such that all trees in the set have their number of plane representations lying in the interval $[e^x,e^{x+C}]$, where $C:=\frac{\log 4}{c}$. Consequently, the number of Pólya trees in the interval must be at least
\begin{equation*}
    K := K(n) = \frac{4^{n-O(\log{n})}}{e^{x+C}} = e^{\log{(4)} n - x - C - O(\log{n})} .
\end{equation*}

If we attach all of these Pólya trees to a common root, the tree we obtain in this way will have $N:=N(n)=Kn+1$ vertices and at least
\begin{align*}
    K! e^{Kx} &= \exp{\left(K\log{K} + Kx - O(K)\right)} \\
    &= \exp{\left(K(\log{(4)} n - x - C - O(\log{n})) + Kx - O(K)\right)}\\
    &=  \exp{\left(\log{(4)}Kn  - O(K\log{n})\right)} 
    = \exp{\left(\log{(4)}N - O\left(\frac{N\log{n}}{n}\right)\right)} 
\end{align*}
plane representations, by Stirling's approximation.

Thus, the probability of picking a plane tree of size $N$ belonging to this isomorphism class when we pick one uniformly at random is
\begin{equation*}
    \frac{\exp{\left(\log{(4)}N - O\left(\frac{N\log{n}}{n}\right)\right)}}{4^{N-O(\log{N})}} = \exp{\left(-O\left(\frac{N(n)\log{n}}{n}\right)\right)} ,
\end{equation*}
which decays subexponentially with $N$. The probability that two trees are isomorphic is at least as large as the probability that they both belong to this isomorphism class, which gives a lower bound of
\begin{equation*}
    \exp{\left(- O\left(\frac{N(n)\log{n}}{n}\right)\right)}^2 = \exp{\left(-O\left(\frac{N(n)\log{n}}{n}\right)\right)} ,
\end{equation*}
for the probability that two plane trees belong to the same isomorphism class which, once again, decays subexponentially. 

\section{Vertex degrees of isomorphic labeled trees}
\label{sec:clt}

We prove Theorem \ref{thm:cltlabeled} with the ultimate goal of comparing the results on the degree distribution of trees conditioned on being isomorphic to the case of standard labeled trees. This is relevant since the conditioning will change the probability distribution. 

Considering two random trees $\mathcal{T}_n^{(1)},\mathcal{T}_n^{(2)}$, we let $Y_{d}^{(i)}$ be a random vector that counts the total number of vertices of degree $\mathbf{d} = (d_1,d_2,\ldots,d_k)$ in tree $i$. For $\mathbf{m} = (m_1,m_2,\ldots,m_k)$ with $0\leq m_i \leq n$, we want to study
\begin{multline*}
    P\left(Y_\mathbf{d}^{(1)} = Y_\mathbf{d}^{(2)} = \mathbf{m}\big|\mathcal{T}_n^{(1)}\simeq \mathcal{T}_n^{(2)}\right) = \frac{P(\{Y_\mathbf{d}^{(1)}=Y_\mathbf{d}^{(2)}=\mathbf{m}\}\cap\{\mathcal{T}_n^{(1)}\simeq\mathcal{T}_n^{(2)}\} )}{P(\mathcal{T}_n^{(1)}\simeq\mathcal{T}_n^{(2)})} \\
    = \frac{\textrm{Total weight of pairs of isomorphic trees with }\mathbf{m}\textrm{ vertices of degree } \mathbf{d}}{\textrm{Total weight of pairs of isomorphic trees}} .
\end{multline*}
We can study this fraction by modifying the functional equations derived in previous sections in a suitable way and then employing general methods from singularity analysis.

Starting from the generating function for pairs of isomorphic labeled trees $n!^2P(x,2)$, we now introduce more variables $\mathbf{u}=(u_1,u_2,\ldots,u_k)$ to keep track of the number of vertices of degree $\mathbf{d}$. The number of vertices of degree $d_i$ is given as the sum of the number of vertices of that degree in each of the root branches plus one if the root is itself of degree $d_i$. Then a modification of the representation (\ref{eq:labeled_expand}), where we recall that we are summing over the root degree, gives the following functional equation for this new multivariate generating function $P(x,2,\mathbf{u})$:
\begin{multline}
\label{eq:MvarFuncEq}
    P(x,2,\mathbf{u})= x\exp\left(\sum_{j=1}^\infty \frac{c(j,2)}{j} P(x^j,2j,\mathbf{u}^j)\right) \\
    + (u_1-1)x\sum_{\lambda \vdash d_1} \prod_j \frac{\left(c(j,2)P(x^j,2j,\mathbf{u}^j)\right)^{\lambda_j}}{j^ {\lambda_j}\lambda_j!} + \ldots \\
    + (u_k-1)x\sum_{\lambda \vdash d_k} \prod_j \frac{\left(c(j,2)P(x^j,2j,\mathbf{u}^j)\right)^{\lambda_j}}{j^ {\lambda_j}\lambda_j!} ,
\end{multline}
where $\mathbf{u}^j=(u_1^j,u_2^j,\ldots,u_k^j)$.
The equation is very similar to the one we had for labeled trees in previous sections; the difference is the extra $k$ terms, which are polynomial in the factors $P(x^j,2j,\mathbf{u}^j)$ for $1\leq j \leq d_k$, where $d_k$ is the largest outdegree we are keeping track of. We then see that the probability generating function of the random vector $\mathbf{X}_{\mathbf{d},n}$ that counts the number of vertices of degree $\mathbf{d}$ in isomorphic pairs of labeled trees, is given by the expression
\begin{equation*}
    \frac{n!^2[x^n]P(x,2,\mathbf{u})}{n!^2[x^n]P(x,2,\mathbf{1})} = \frac{[x^n]P(x,2,\mathbf{u})}{[x^n] P(x,2)} .
\end{equation*}

As the functions with $j\geq 2$ are slight modifications of analytic functions, we have reason to think that equation (\ref{eq:MvarFuncEq}) should be well behaved. We have already seen in Section \ref{sec:labeled} that $\alpha$, the singularity of $P(x,2)$, satisfies $0<\alpha<1$, so we find that $0<|\alpha u_1u_2\cdots u_k|<1$ if we restrict $\mathbf{u}$ to some ball around $\mathbf{1}$, i.e. if we assume that $|\mathbf{u}-\mathbf{1}|<\epsilon_1$ for some $\epsilon_1>0$. Since we have $0\leq d_i\leq n$ for all $i$, we can then argue in a similar way as we did in Section \ref{sec:labeled} that $P(x^j,2j,\mathbf{u}^j)$ is analytic in a larger region than $P(x^j,2,\mathbf{u}^j)$. We can thereby conclude that the higher order factors $P(x^j,2j,\mathbf{u}^j)$ for $j>1$ are analytic in a region $|x|<\alpha+\epsilon_2$ for some $\epsilon_2>0$. This lets us apply the following theorem.
\begin{theorem}[Theorem 2.23 and Remark 2.24 from \cite{DrmotaBook_MR2484382}]
\label{thm:combCLT}
Let $\mathbf{u} = (u_1,u_2,\ldots,u_k)$ and suppose that $\mathbf{X}_n$ is a sequence of random vectors such that
\begin{equation*}
    \E \mathbf{u}^{\mathbf{X}_n} = \frac{[x^n]y(x,\mathbf{u})}{[x^n]y(x,\mathbf{1})} ,
\end{equation*}
where $y(x, \mathbf{u})$ is a power series, that is the (analytic) solution of the functional equation $y = F(x, y, \mathbf{u})$, where $F(x, y, \mathbf{u}) = \sum_{n,m} F_{n,m}(\mathbf{u})$ is an analytic function in $x, y$ around $0$ and $\mathbf{u}$ around $\mathbf{1}$ such that $F(0, y, u) \equiv 0$, that
$F(x, 0, u) \equiv 0$, and that all coefficients $F_{n,m}(\mathbf{1})$ of $F(x, y, \mathbf{1})$ are real and non-negative. 

Assume that the region of convergence of $F(x,y,\mathbf{u})$ is large enough for there to exist a non-negative solution $x = x_0 > 0$ and $y = y_0 > 0$ to the system of equations
\begin{align*}
    y = F(x,y,\mathbf{1}) ,\\
    1 = F_y(x,y,\mathbf{1}) ,
\end{align*}
with $F_x(x_0,y_0,\mathbf{1}) \neq 0$ and $F_{yy}(x_0,y_0,\mathbf{1})\neq 0$. Then we get
\begin{equation*}
    \E \mathbf{X}_n = \boldsymbol{\mu} n + O(1)\textrm{ and } \cov\mathbf{X}_n = \mathbf{\Sigma} n + O(1),  
\end{equation*}
where $\boldsymbol{\mu} = (\mu_1, \mu_2,\ldots,\mu_k)$ and $\mathbf{\Sigma} = (\sigma_{i,j})_{1\leq i,j\leq k}$ can be calculated. E.g., for the mean we have
\begin{equation*}
    \mu_i = \frac{F_{u_i}}{x_0F_x} ,
\end{equation*}
where all partial derivatives are evaluated at the point $(x_0, y_0, \mathbf{1})$. Furthermore, we have a central limit theorem of the form
\begin{equation*}
    \frac{1}{\sqrt{n}}\left(\mathbf{X}_n - \E \mathbf{X}_n\right) \xrightarrow{d} N(\mathbf{0},\mathbf{\Sigma}) .
\end{equation*}
\end{theorem}

We can now write the equation (\ref{eq:MvarFuncEq}) as $y = F(x,y,\mathbf{u})$ where we have exchanged $P(x,2,\mathbf{u})$ with $y$. Then, the second part is a polynomial in analytic functions and the variables $x,y,\mathbf{u}$ and thus analytic itself, and, by similar reasoning as in Section \ref{sec:labeled}, the first part is analytic in a region containing $x=x_0 = \alpha$, $y=y_0=P(\alpha,2,\mathbf{1})$ as long as we restrict $\mathbf{u}$ to be close to $\mathbf{1}$. This means that the functional equation satisfies the conditions in Theorem \ref{thm:combCLT} and we find that the random variable $\mathbf{X}_{\mathbf{d},n}$ asymptotically follows a multivariate normal distribution.

\subsection{Example: leaves in isomorphic labeled trees.}
\label{subsec:Leaveslabeled}
For the special case of leaves, i.e., when $X_n$ is a random variable that counts vertices with (out-)degree 0, the functional equation becomes
\begin{equation*}
    P(x,2,u) = x\exp\left(\sum_{j=1}^\infty \frac{c(j,2)}{j} P(x^j,2j,u^j)\right) + (u-1)x ,
\end{equation*}
and we can calculate the mean as
\begin{equation*}
     \E X_{n} = \frac{[x^n]P_u(x,2,1)}{[x^n]P(x,2,1)} .
\end{equation*}
By differentiating the equation above, we find
\begin{equation*}
    P_u(x,2,1) = \frac{1}{1-P(x,2,1)} \left(P(x,2,1)\sum_{j\geq 2} c(j,2) P_u(x^j,2j,1) + x \right) .
\end{equation*}
Now set 
\begin{equation*}
    A(x) = \sum_{j\geq 2} c(j,2) P_u(x^j,2j,1) ,
\end{equation*}
then, after performing singularity analysis and combining this with the asymptotics we derived for $[x^n]P(x,2,1)=[x^n]P(x,2)$ in Section \ref{sec:labeled} we find that
\begin{equation*}
    \E X_{n} = \mu n + O(1),
\end{equation*}
where, using the same notation as in Section \ref{sec:labeled}, 
\begin{equation*}
    \mu = \left(\frac{A(\alpha)}{\alpha\xi'(\alpha)} + \frac{1}{\xi'(\alpha)}\right) e^{-1} .
\end{equation*} 
Proceeding as in Section \ref{sec:labeled}, we can estimate $\mu\approx 0.340252$. This can be compared to the case of ordinary labeled trees, where the mean constant is $e^{-1} \approx 0.367879$ (see \cite[Theorem 3.13]{DrmotaBook_MR2484382}). Thus, a pair of labeled trees conditioned on being isomorphic exhibits a different shape than the usual model of a uniformly random labeled tree.


\section{Weights of isomorphism classes of Pólya trees with restricted degrees}
\label{subsec:nrRep}

The functional equation (\ref{eq:funcEqBoundDeg}) describes the weights of isomorphism classes of trees with bounded degrees. Recall that the variable $t$ in $P_D(x,t)$ keeps track of $w(T)\mathcal{PR}(P) = W(P) = e^{\log{W(P)}}$, where $T$ is a simple tree belonging to the isomorphism class corresponding to the Pólya tree $P$. We want to apply a generalization of Theorem \ref{thm:combCLT}, where we also make the substitution $u\to e^t$, to this equation to derive a central limit theorem for $\log{W(\mathcal{P}_n)}$, i.e. to prove Theorem \ref{thm:nrRepGW}. To be precise, we generalize Theorem \ref{thm:combCLT} to allow the additive parameters, i.e. the weights of the isomorphism classes in our case, to be non-integers. The proof of Theorem \ref{thm:combCLT} is general enough to allow for this modification. To verify that the theorem really applies here, we need to check that the functions $P_D(x^j,jt)$ for $j\geq 2$ involved in the expression are analytic in some region containing $\rho_D$, the singularity of $P_D(x,0)$, as long as we restrict $t$ to lie in some small disc around 0. As the right hand side of equation (\ref{eq:funcEqBoundDeg}) is a polynomial in these functions and $P(x,t)$, this would imply that the expression is analytic in a suitably large region for the method of singularity analysis to apply. 

Let $\rho_D$ be the singularity of the Pólya trees corresponding to the isomorphism classes of the trees $\mathcal{T}$ in question. This number is known to be less than 1 (except if $D=\{0,1\}$), since the number of trees of size $n$ grows exponentially. Indeed, we can construct exponentially many different trees by starting with a rooted path consisting of $m$ vertices of degree $d$, where $m$ is an integer to be decided and $d>0$ is some allowed outdegree in the tree. We then pick out $\frac{m}{2}$ vertices on the path and, for each of them, we pick out one of its children and attach $d$ vertices to it. There are 
\begin{equation*}
\binom{m}{\frac{m}{2}} \sim C m^{-\frac{1}{2}} 2^m 
\end{equation*}
possible choices and each of them results in a unique tree. Setting $m = \frac{2(n-1)}{3d}$ forces the total size of the tree to be $n$ and shows that the total number of trees of size $n$ grows (at least) exponentially in $n$. That the number of trees does not grow faster than exponential follows from the fact that they cannot grow faster than the total number of Pólya trees (without degree restrictions).

We want to check that $|P_D(x^j,jt)|$ is finite in a region $|x|<\rho_D + \epsilon$, for some $\epsilon>0$, as long as we restrict $t$ to be close to 0. We know that $\rho_D<1$, so we only need to check that the factors $|W(P)^t|$ do not grow too fast with the size of $P$. As we have
\begin{equation*}
    |W(P)^{a+bi}| = W(P)^a
\end{equation*}
it is enough to consider the cases where $t$ is real. First assume that $t$ is positive. Recall that the weight of simply generated trees (under some fairly mild conditions) behaves like
\begin{equation*}
\label{eq:expDecayGW}
    [x^n]T(x) \sim C n^{-3/2} \rho^{-n} 
\end{equation*}
asymptotically. As a consequence, we have a bound of the same type for the isomorphism classes of trees -- no isomorphism class can have larger weight than the total sum of weights. Therefore, we can find some $K$ such that
\begin{equation*}
    W(P)^{t} \leq K^{tn} .
\end{equation*}
 If we instead assume that $t=-a$ is negative we can use the fact that the set $D$ is finite to see that
\begin{equation*}
    \frac{1}{W(P)^a} = \frac{1}{w(T)^a \mathcal{PR}(P)^a} \leq \prod_v w_{\deg(v)}^{-a} \leq \min_{d\in D}\{w_d\}^{-an} = \min_{d\in D}\{w_d\}^{tn} .
\end{equation*}
In either case, $|W(P)|^t$ is at most exponential in $|P|$, and, by taking $t$ to be close enough to $0$, we can guarantee that $P_D(x^j,jt)$, where $j$ is bounded by the largest out-degree, will be analytic in a region containing $\rho_D$. Since this holds, Theorem \ref{thm:combCLT} applies directly and we get the central limit theorem with both mean and variance being of order $n$. 

We can also prove Theorem \ref{thm:nrRepLab}, i.e. a central limit theorem for isomorphism classes of labeled trees (without any degree restrictions), more or less for free, based on the results in \cite{autRandTrees}. 
In this context, the weight of the isomorphism classes corresponds to the number of labelings of the underlying Pólya tree. The number of labelings of a Pólya tree is $\frac{n!}{|\aut{P}|}$ so when we take the logarithm we get
\begin{equation*}
    \log n! - \log|\aut{P}| .
\end{equation*}
The first term is deterministic and the second satisfies a central limit theorem, as was shown in \cite{autRandTrees}. Using Stirling's approximation, we find that the number of labelings of a random Pólya tree asymptotically follows a normal distribution with expected value
\begin{equation*}
    \log n!-\E[\log|\aut{P}|] = n\log n - (\mu+1) n + \frac{1}{2} \log n + O(1) ,
\end{equation*}
for some constant $\mu \approx 0.137342$ and variance $\sigma^2 n + O(1)$, where $\sigma^2 \approx 0.196770$. Note how the order of the mean is $n\log n$ which differs from what we had for Galton--Watson trees with bounded degrees (as well as many other examples in combinatorial probability due to the strength of Theorem \ref{thm:combCLT}).

\subsection{Example: binary trees} 
\label{subsec:EstPlaneRepBin} 
As a special case we consider the relationship between unordered, unlabeled binary trees with $D=\{0,1,2\}$, and binary simple trees with weights $\mathbf{w} = [1,2,1]$. These classes of trees satisfy Theorem \ref{thm:nrRepGW}, which means that the number of plane representations of these trees, where we always distinguish between left and right children in the embedding, is asymptotically log-normal. The functional equation (\ref{eq:funcEqBoundDeg}) becomes 
\begin{equation*}
    P_D(x,t) = x + 2^txP_D(x,t) + 2^{t-1} xP_D(x,t)^2 + (1-2^{t-1})x P_D(x^2,2t)
\end{equation*}
in this case. If we truncate the rapidly converging function $P_{\{0,1,2\}}(x^2,2t)$ after a few terms, we can estimate the constants numerically with general formulas provided by Theorem \ref{thm:combCLT} (slightly modified to take the change of variables $u\to e^t$ into account). We find that the expected value is
\begin{equation*}
    \E[\log{W(\mathcal{P}_n)}] = \mu n + O(1) ,
\end{equation*}
where $\mu \approx 0.444518$, and that the variance is
\begin{equation*}
    \Var[\log{W(\mathcal{P}_n)}] = \sigma^2 n + O(1) ,
\end{equation*}
where $\sigma^2\approx 0.072413$.

If we instead do the same where the Pólya trees are embedded as unary-binary trees (with weights $\mathbf{w} = [1,1,1]$), so that we do not differentiate between right and left children for vertices that only have one child, we again have a central limit theorem for $\log W(\mathcal{P}_n)$. In this case we have $\mu \approx 0.176278$ and $\sigma^2 \approx 0.025865$. In other words, the average number of embeddings is smaller, just as one would expect. 


\section{Acknowledgments}
The author would like to thank Stephan Wagner for many helpful discussions and suggestions during the writing of this article.

\bibliographystyle{unsrt}
\bibliography{references}

\end{document}